\documentclass{amsart}
\usepackage{setspace}
\usepackage{a4}
\usepackage{amsthm}
\usepackage{latexsym}
\usepackage{amsfonts}
\usepackage{graphicx}
\usepackage{textcomp}
\usepackage{cite}
\usepackage{enumerate}
\usepackage{amssymb}
\usepackage{hyperref}
\usepackage{amsmath}
\usepackage{tikz}
\usepackage[mathscr]{euscript}
\usepackage{mathtools}
\newtheorem{theorem}{Theorem}[section]

\newtheorem{corollary}[theorem] {Corollary}
\newtheorem{definition}[theorem]{Definition}

\newtheorem{lemma} [theorem]{Lemma}

\newtheorem{problem}[theorem]{Problem}

\newtheorem{remark}[theorem]{Remark}

\title{This is the title}
\usepackage{fancyhdr}

\begin{document}
\hrule\hrule\hrule\hrule\hrule
\vspace{0.3cm}	
\begin{center}
{\bf{A Functional Version of the Sparsity Theorem}}\\
\vspace{0.3cm}
\hrule\hrule\hrule\hrule\hrule
\vspace{0.3cm}
\textbf{K. Mahesh Krishna}\\
School of Mathematics and Natural Sciences\\
Chanakya University Global Campus\\
NH-648, Haraluru Village\\
Devanahalli Taluk, 	Bengaluru  North District\\
Karnataka State 562 110 India\\
Email: kmaheshak@gmail.com\\

Date: \today
\end{center}

\hrule\hrule
\vspace{0.5cm}
\textbf{Abstract}: Celebrated breakthrough sparsity theorem obtained independently by Donoho and Elad \textit{[Proc. Natl. Acad. Sci. USA, 2003]} and Gribonval and Nielsen \textit{[IEEE Trans. Inform. Theory, 2003]} and Fuchs \textit{[IEEE Trans. Inform. Theory, 2004]} says that unique  sparse solution to NP-Hard $\ell_0$-minimization problem can be obtained using unique solution to P-Type $\ell_1$-minimization problem. In this paper, we extend their result to abstract Banach spaces.  We notice that the `normalized' condition for Hilbert spaces can be generalized to a larger extent when we consider Banach spaces.

\textbf{Keywords}:   Compressed sensing, Frame, Schauder frame, sparse solution.

\textbf{Mathematics Subject Classification (2020)}: 42C15, 94A12, 46B45.\\

\hrule

\hrule
\section{Introduction}

Let $\mathcal{H}$ be a  finite dimensional Hilbert  space over $\mathbb{K}$ ($\mathbb{C}$ or $\mathbb{R}$). Recall  \cite{BEDETTOFICKUS, HANKORNELSON} that a finite collection $\{\tau_j\}_{j=1}^n$ in $\mathcal{H}$ is said to be a \textbf{frame} (also known as \textbf{dictionary}) for $\mathcal{H}$ if it spans $\mathcal{H}$. A frame $\{\tau_j\}_{j=1}^n$ for  $\mathcal{H}$ is said to be \textbf{normalized} if $\|\tau_j\|=1$ for all $1\leq j \leq n$. Given a frame $\{\tau_j\}_{j=1}^n$ for  $\mathcal{H}$, we define the analysis operator
\begin{align*}
	\theta_\tau:\mathcal{H}\ni h \mapsto 	\theta_\tau h \coloneqq(\langle h,  \tau_j \rangle)_{j=1}^n \in \mathbb{K}^n.
\end{align*}
Adjoint of the analysis operator is known as  the synthesis operator whose equation is 
\begin{align*}
\theta_\tau^*: \mathbb{K}^n \ni (a_j)_{j=1}^n \mapsto \theta_\tau^*(a_j)_{j=1}^n\coloneqq\sum_{j=1}^{n}a_j\tau_j \in \mathcal{H}.
\end{align*}
Given $d\in \mathbb{K}^n$, let $\|d\|_0$ be the number of nonzero entries in $d$. Problems that arise in single-pixel cameras, magnetic resonance imaging (MRI), radar, speech recognition, coding and information theory, machine learning, hyper-spectral imaging, geophysical data analysis, remote sensing, robotics and control, analog-to-digital conversion, astronomy etc can be formulated as the following $\ell_0$-minimization problem.
\begin{problem} \label{P0}
 Let $\{\tau_j\}_{j=1}^n$ be a frame for  $\mathcal{H}$.  Given $h \in \mathcal{H}$, solve 	
\begin{align*}
	\underset{d\in \mathbb{K}^n}{\operatorname{minimize}}\,\|d\|_0\quad \text{ subject to } \quad \theta_\tau^*d=h.
\end{align*}
\end{problem}
Recall \cite{BEDETTOFICKUS, HANKORNELSON} that $c \in \mathbb{K}^n$ is said to be a unique solution to Problem \ref{P0}  if it satisfies the following two conditions.
\begin{enumerate}[\upshape(i)]
	\item $\theta_\tau^*c=h$.
	\item If $d \in \mathbb{K}^n$ satisfies $\theta_\tau^*d=h$, then 
\begin{align*}
	\|d\|_0>\|c\|_0.
\end{align*}
\end{enumerate}
Unfortunately, in 1995, Natarajan showed that Problem \ref{P0} is NP-Hard \cite{NATARAJAN, FOUCARTRAUHUT}. Therefore the solution to Problem \ref{P0} has to be obtained using other means, i.e., using solution to P-problem which is solvable using linear programmings. Entire body of work which is built around Problem \ref{P0} is known as \textbf{sparseland}  (term due to Elad \cite{ELAD}) or \textbf{compressive sensing} or \textbf{compressed sensing} \cite{VIDYASAGAR, ELAD, DAVENPORTDUARTEELDARKUTYNIOK, FOUCARTRAUHUT, KUTYNIOK2, DONOHO, ADCOCKHANSEN, RISHGRABARNIK, ELDARBOOK, TILLMANNPFETSCH, COHENDAHMENDEVORE, CANDESTAO2, CANDESTAO, CANDESTAO3, CANDESTAO4, CANDESROMBERGTAO, CANDES, CANDESROMBERGTAO2, FEUERNEMIROVSKI, GRIBONVALNIELSEN, DONOHOHUO, DONOHOELAD, ELADBRUCKSTEIN}.
We  note that as the operator $\theta_\tau^*$ is surjective, for a given $h \in \mathcal{H}$, there is always $d\in \mathbb{K}^n$ such that  $\theta_\tau^*d=h.$ Thus the central problem is when the solution to Problem \ref{P0} is unique.   It is well-known that (see \cite{CHENDONOHOSAUNDERS, DONOHOHUO, BRUCKSTEINDONOHOELAD}) the  following problem is the closest convex relaxation problem to Problem \ref{P0}.
\begin{problem} \label{P1}
Let $\{\tau_j\}_{j=1}^n$ be a frame for  $\mathcal{H}$.	Given $h \in \mathcal{H}$, solve 	
	\begin{align*}
		\underset{d\in \mathbb{K}^n}{\operatorname{minimize}}\,\|d\|_1 \quad \text{ subject to } \quad \theta_\tau^*d=h.
	\end{align*}
\end{problem}
Problem 1.2 is a P-problem, but several solutions are available via linear programming \cite{XUEYE, TERLAKY, TILLMANN}.

The following important result shows that the solution to
Problem \ref{P1} gives the solution to Problem \ref{P0}. It was obtained independently by Donoho and Elad \cite{DONOHOELAD} and Gribonval and Nielsen \cite{GRIBONVALNIELSEN} and Fuchs \cite{FUCHS, FUCHS2}.
\begin{theorem}\cite{DONOHOELAD, GRIBONVALNIELSEN, ELAD, KUTYNIOK, FUCHS, FUCHS2}\label{DEGN} (\textbf{Donoho-Elad-Gribonval-Nielsen-Fuchs Sparsity Theorem})
	Let $\{\tau_j\}_{j=1}^n$ be a  normalized frame  for  a Hilbert space $\mathcal{H}$. If  
 $h \in \mathcal{H}$ can be written as  $	h=\theta_\tau^*c$ for some  $c\in \mathbb{K}^n$ satisfying 
 
 \begin{align*}
 \|c\|_0<\frac{1}{2}\left(1+\frac{1}{\displaystyle\max_{1\leq j< k \leq n}|\langle \tau_j, \tau_k\rangle|}\right),
 \end{align*}
 then $c$ is the unique solution to Problem \ref{P1} and Problem \ref{P0}.
\end{theorem}
We naturally ask for (both finite and infinite dimensional)  Banach space version  of Theorem \ref{DEGN}.  As frame theory for Hilbert spaces has been successfully extended to Banach spaces, Theorem  \ref{DEGN} has not yet found a generalization. It is interesting to note that a noncommutative version of Theorem \ref{DEGN} (using frames for Hilbert C*-modules) has been recently derived \cite{KRISHNA}.

\section{Functional Donoho-Elad-Gribonval-Nielsen-Fuchs Sparsity Theorem}
In the paper,   $\mathcal{X}$ denotes a  Banach space (need not be finite dimensional) over $\mathbb{K}$. Dual of $\mathcal{X}$ is denoted by $\mathcal{X}^*$.
We begin by formulating  Problems \ref{P0} and \ref{P1} in general forms. 
\begin{problem} \label{BP0}
Let $\{\tau_n \}_{n=1}^\infty$ be a sequence in  $\mathcal{X}$ such that the map 
\begin{align*}
	\theta_\tau:\ell^1(\mathbb{N})\ni \{a_n \}_{n=1}^\infty\mapsto \theta_\tau \{a_n \}_{n=1}^\infty\coloneqq \sum_{n=1}^{\infty}a_n\tau_n \in \mathcal{X} 
\end{align*}
is a well-defined bounded linear surjective operator. Given $x \in \mathcal{X}$, solve 	
	\begin{align*}
	\underset{d\in \ell^1(\mathbb{N})}{\operatorname{minimize}}\,\|d\|_0 \quad \text{ subject to } \quad \theta_\tau d=x.
\end{align*}
\end{problem}
\begin{problem} \label{BP1}
	Let $\{\tau_n \}_{n=1}^\infty$ be a sequence in  $\mathcal{X}$ such that the map 
	\begin{align*}
		\theta_\tau:\ell^1(\mathbb{N})\ni \{a_n \}_{n=1}^\infty\mapsto \theta_\tau \{a_n \}_{n=1}^\infty\coloneqq \sum_{n=1}^{\infty}a_n\tau_n \in \mathcal{X} 
	\end{align*}
	is a well-defined bounded linear surjective operator. Given $x \in \mathcal{X}$, solve 	
	\begin{align*}
		\underset{d\in \ell^1(\mathbb{N})}{\operatorname{minimize}}\,\|d\|_1 \quad \text{ subject to } \quad \theta_\tau d=x.
	\end{align*}
\end{problem}
A very important property used to show Theorem \ref{DEGN} is the null space property (see \cite{KUTYNIOK, COHENDAHMENDEVORE}). We now define an analogous property for Banach spaces. We use the following notations. Let $\{e_n \}_{n=1}^\infty$ be the canonical Schauder basis for $\ell^1(\mathbb{N})$. Given $M\subseteq \mathbb{N}$ and $d=\{d_n \}_{n=1}^\infty \in \ell^1(\mathbb{N})$, we define 
\begin{align*}
	d_M\coloneqq \sum_{n\in M}d_ne_n.
\end{align*}
Number of elements in  $ M \subseteq \mathbb{N}$ is denoted by   $|M|$.
\begin{definition}
Let $\{\tau_n \}_{n=1}^\infty$ be a sequence in  $\mathcal{X}$ such that the map 
\begin{align*}
	\theta_\tau:\ell^1(\mathbb{N})\ni \{a_n \}_{n=1}^\infty\mapsto \theta_\tau \{a_n \}_{n=1}^\infty\coloneqq \sum_{n=1}^{\infty}a_n\tau_n \in \mathcal{X} 
\end{align*}
is a well-defined bounded linear surjective operator. The sequence $\{\tau_n \}_{n=1}^\infty$  is said to have the \textbf{null space property} (we write NSP) of order $k\in \mathbb{N}$ if 	for every $ M \subseteq \mathbb{N} $ with $|M|\leq k$, we have
\begin{align*}
	\|d_M\|_1<\frac{1}{2}\|d\|_1,  \quad  \forall d \in \ker(\theta_\tau), d \neq 0.
\end{align*}
\end{definition}
The following characterization relates NSP with Problem \ref{BP1}.
\begin{lemma}\label{IFF}
Let $\{\tau_n \}_{n=1}^\infty$ be a sequence in  $\mathcal{X}$ such that the map 
\begin{align*}
	\theta_\tau:\ell^1(\mathbb{N})\ni \{a_n \}_{n=1}^\infty\mapsto \theta_\tau \{a_n \}_{n=1}^\infty\coloneqq \sum_{n=1}^{\infty}a_n\tau_n \in \mathcal{X} 
\end{align*}
is a well-defined bounded linear surjective operator. Let $k \in \mathbb{N}$.
 The following are equivalent.
\begin{enumerate}[\upshape(i)]
	\item  If $x \in \mathcal{X}$ can be written as  $	x=\theta_\tau c$ for some  $c \in \ell^1(\mathbb{N})$ with $\|c\|_0\leq k$, then $c$ is the unique solution to Problem \ref{BP1}.
	\item $ \{\tau_n \}_{n=1}^\infty$  satisfies the  null space property  of order $k$.
\end{enumerate}	
\end{lemma}
\begin{proof}
\begin{enumerate}[\upshape(i)]
	\item 	$\implies$ (ii) Let $ M \subseteq \mathbb{N} $ with $|M|\leq k$ and let $d \in \ker(\theta_\tau), d \neq 0$. Then we have 
	\begin{align*}
	0=\theta_\tau d=\theta_\tau(d_M+d_{M^c})=\theta_\tau(d_M)+\theta_\tau(d_{M^c})
	\end{align*}
which gives 
\begin{align*}
	\theta_\tau(d_M)=\theta_\tau(-d_{M^c}).
\end{align*}
Define $c\coloneqq d_M \in \ell^1(\mathbb{N})$ and $x\coloneqq \theta_\tau(d_M)$. Then we have $\|c\|_0\leq |M|\leq k$ and 
\begin{align*}
	x=\theta_\tau c=\theta_\tau(-d_{M^c}).
\end{align*}
By assumption (i),  we then have 
\begin{align*}
	\|c\|_1=\|d_M\|_1<\|-d_{M^c}\|_1=\|d_{M^c}\|_1.
\end{align*}
Rewriting previous inequality gives 
\begin{align*}
	\|d_M\|_1<\|d\|_1-\|d_M\|_1 \implies 	\|d_M\|_1<\frac{1}{2}\|d\|_1.
\end{align*}
Hence $\{\tau_n \}_{n=1}^\infty$  satisfies the  NSP of order $k$.
	\item    $\implies$ (i)  Let  $x \in \mathcal{X}$ can be written as  $	x=\theta_\tau c$ for some  $c \in \ell^1(\mathbb{N})$ satisfying $\|c\|_0\leq k$. Define $M\coloneqq \operatorname{supp}(c)$. Then $|M|=\|c\|_0 \leq k$. By assumption (ii), we then have 
	\begin{align}\label{PI}
	\|d_M\|_1<\frac{1}{2}\|d\|_1,  \quad  \forall d \in \ker(\theta_\tau), d \neq 0.
\end{align}
Let $b\in \ell^1(\mathbb{N})$ be such that $x=\theta_\tau b$ and $b\neq c$. Define $a\coloneqq b-c \in \ell^1(\mathbb{N})$. Then $\theta_\tau a=\theta_\tau b-\theta_\tau c=x-x=0$  and hence $a \in \ker(\theta_\tau), a \neq 0$. Using inequality (\ref{PI}), we get 
\begin{align}\label{P2}
	\|a_M\|_1<\frac{1}{2}\|a\|_1 \implies \|a_M\|_1<\frac{1}{2}(\|a_{M}\|_1+\|a_{M^c}\|_1) \implies \|a_M\|_1<\|a_{M^c}\|_1.
\end{align}
Using inequality (\ref{P2}) and the information that $c$ is supported on $M$, we get 
\begin{align*}
\|b\|_1-\|c\|_1&=\|b_M\|_1+\|b_{M^c}\|_1-\|c_M\|_1-\|c_{M^c}\|_1=\|b_M\|_1+\|b_{M^c}\|_1-\|c_M\|_1\\
&=\|b_M\|_1+\|(b-c)_{M^c}\|_1-\|c_M\|_1=\|b_M\|_1+\|a_{M^c}\|_1-\|c_M\|_1\\
&>\|b_M\|_1+\|a_M\|_1-\|c_M\|_1\geq \|b_M\|_1+\|(b-c)_M\|_1-\|c_M\|_1\\
&\geq \|b_M\|_1-\|b_M\|_1+\|c_M\|_1-\|c_M\|_1= 0.
\end{align*}
Hence $c$ is the unique solution to Problem \ref{BP1}.
\end{enumerate}	
\end{proof}
Using Theorem \ref{IFF} we obtain Banach space version of Theorem \ref{DEGN}. We do this  by relating Problem \ref{BP1} to Lemma \ref{IFF} and then Problem \ref{BP0} to Lemma \ref{IFF}. Note that the following theorem requires the existence of a new sequence of functionals.
\begin{theorem}\label{TDP}
	Let $\{\tau_n \}_{n=1}^\infty$ be a sequence in  $\mathcal{X}$ such that the map 
	\begin{align*}
		\theta_\tau:\ell^1(\mathbb{N})\ni \{a_n \}_{n=1}^\infty\mapsto \theta_\tau \{a_n \}_{n=1}^\infty\coloneqq \sum_{n=1}^{\infty}a_n\tau_n \in \mathcal{X} 
	\end{align*}
	is a well-defined bounded linear surjective operator. Assume that there is a sequence $\{f_n \}_{n =1}^\infty$ in $\mathcal{X}^*$ such that the map 
	\begin{align}\label{L}
		\theta_f: \mathcal{X} \ni x \mapsto \theta_fx\coloneqq \{f_n(x) \}_{n =1}^\infty \in \ell^1(\mathbb{N})
	\end{align}
	is a well-defined bounded linear operator and 
	 \begin{align}\label{N}
	|f_n(\tau_n)|\geq1, \quad \forall n \in \mathbb{N}.
\end{align}
	If $x \in \mathcal{X}$ can be written as  $	x=\theta_\tau c$ for some  $c \in \ell^1(\mathbb{N})$ satisfying 
\begin{align}\label{I}
\|c\|_0<\frac{1}{2}\left(1+\frac{1}{\displaystyle\sup_{n, m \in \mathbb{N}, n \neq m}| f_n(\tau_m)|}\right),	
\end{align}
then $c$ is the unique solution to Problem \ref{BP1}.
\end{theorem}
\begin{proof}
We show that $\{\tau_n \}_{n=1}^\infty$ satisfies the  NSP of order $k\coloneqq \|c\|_0$. Then Lemma \ref{IFF} says that $c$ is the unique solution to Problem \ref{BP1}. Let  $x \in \mathcal{X}$  be written as  $	x=\theta_\tau c$ for some  $c \in \ell^1(\mathbb{N})$ satisfying $\|c\|_0\leq k$. Let $ M \subseteq \mathbb{N} $ with $|M|\leq k$ and let $d \in \ker(\theta_\tau), d \neq 0$. Then we have 
\begin{align*}
	\theta_f\theta_\tau d=0.
\end{align*}
By writing  $d=\{d_n\}_{n=1}^\infty\in \ell^1(\mathbb{N})$, the above equation gives 
\begin{align*}
0&=\theta_f\theta_\tau\{d_m\}_{m=1}^\infty=\theta_f\left(\sum_{m=1}^{\infty}d_m\theta_\tau e_m\right)\\
&=\theta_f\left(\sum_{m=1}^{\infty}d_m\tau_m\right)=
\sum_{m=1}^{\infty}d_m\theta_f(\tau_m)=\sum_{m=1}^{\infty}d_m\sum_{k=1}^{\infty}f_k(\tau_m)e_k.
\end{align*}
Let $\{\zeta_n \}_{n=1}^\infty$ be the coordinate functionals associated with the canonical Schauder basis  $\{e_n \}_{n=1}^\infty$ for   $\ell^1(\mathbb{N})$. Let $n \in \mathbb{N} $. By evaluating previous equation at $\zeta_n$, we get 

	 \begin{align*}
0&=\zeta_n\left(\sum_{m=1}^{\infty}d_m\sum_{k=1}^{\infty}f_k(\tau_m)e_k\right)=\sum_{m=1}^{\infty}d_m\sum_{k=1}^{\infty}f_k(\tau_m)\zeta_n(e_k)\\
&=\sum_{m=1}^{\infty}d_mf_n(\tau_m)=d_nf_n(\tau_n)+\sum_{m=1, m\neq n}^{\infty}d_mf_n(\tau_m).
\end{align*}
Therefore 
\begin{align*}
d_nf_n(\tau_n)=-\sum_{m=1, m\neq n}^{\infty}d_mf_n(\tau_m), \quad \forall n \in \mathbb{N}.	
\end{align*}
By using inequality (\ref{N}), 
\begin{align*}
	|d_n|&\leq |d_n||f_n(\tau_n)|=\left|-\sum_{m=1, m\neq n}^{\infty}d_mf_n(\tau_m)\right|\\
	&\leq \sum_{m=1, m\neq n}^{\infty}|d_mf_n(\tau_m)|\leq \left(\displaystyle\sup_{m \in \mathbb{N}, n \neq m}|f_n(\tau_m)|\right)\sum_{m=1, m\neq n}^{\infty}|d_m|\\
	&\leq  \left(\displaystyle\sup_{n, m \in \mathbb{N}, n \neq m}|f_n(\tau_m)|\right)\sum_{m=1, m\neq n}^{\infty}|d_m|=\left(\displaystyle\sup_{n, m \in \mathbb{N}, n \neq m}|f_n(\tau_m)|\right)\left(\sum_{m=1}^{\infty}|d_m|-|d_n|\right)\\
	&=\left(\displaystyle\sup_{n, m \in \mathbb{N}, n \neq m}|f_n(\tau_m)|\right)(\|d\|_1-|d_n|), \quad \forall n \in \mathbb{N}.
\end{align*}
By rewriting above inequality  we get
\begin{align}\label{S}
	\left(1+\frac{1}{\displaystyle\sup_{n, m \in \mathbb{N}, n \neq m}| f_n(\tau_m)|}\right)|d_n|\leq \|d\|_1, \quad \forall n \in \mathbb{N}.
\end{align}
Summing inequality (\ref{S}) over $M$ leads to 
\begin{align*}
	\left(1+\frac{1}{\displaystyle\sup_{n, m \in \mathbb{N}, n \neq m}| f_n(\tau_m)|}\right)\|d_M\|_1&=\left(1+\frac{1}{\displaystyle\sup_{n, m \in \mathbb{N}, n \neq m}| f_n(\tau_m)|}\right)\sum_{n\in M}|d_n|\\
	&\leq \|d\|_1	\sum_{n\in M} 1= \|d\|_1|M|.
\end{align*}
Finally using inequality (\ref{I}), we have 
\begin{align*}
\|d_M\|_1&\leq \left(1+\frac{1}{\displaystyle\sup_{n, m \in \mathbb{N}, n \neq m}| f_n(\tau_m)|}\right)^{-1} \|d\|_1|M|\leq\left(1+\frac{1}{\displaystyle\sup_{n, m \in \mathbb{N}, n \neq m}| f_n(\tau_m)|}\right)^{-1} \|d\|_1k\\
&=\left(1+\frac{1}{\displaystyle\sup_{n, m \in \mathbb{N}, n \neq m}| f_n(\tau_m)|}\right)^{-1} \|d\|_1\|c\|_0<\frac{1}{2}\|d\|_1.
\end{align*}
Hence 	$\{\tau_n \}_{n=1}^\infty$ satisfies the  NSP of order $k$ which completes the proof.
\end{proof}
We observe that only the inequality (\ref{N}) was used in the proof of the previous theorem, and equality is not required. This is the precise generality gained when moving from Hilbert spaces to Banach spaces.
\begin{remark}
	We used the existence of functionals satisfying (\ref{L}) and inequality (\ref{N}) in the proof of Theorem \ref{TDP}. We are unable to derive the result without this assumption. 
\end{remark}
\begin{theorem} \label{M}
(\textbf{Functional Donoho-Elad-Gribonval-Nielsen-Fuchs Sparsity Theorem})\\
	Let $\{\tau_n \}_{n=1}^\infty$ be a sequence in  $\mathcal{X}$ such that the map 
\begin{align*}
	\theta_\tau:\ell^1(\mathbb{N})\ni \{a_n \}_{n=1}^\infty\mapsto \theta_\tau \{a_n \}_{n=1}^\infty\coloneqq \sum_{n=1}^{\infty}a_n\tau_n \in \mathcal{X} 
\end{align*}
is a well-defined bounded linear surjective operator. Assume that there is a sequence $\{f_n \}_{n =1}^\infty$ in $\mathcal{X}^*$ such that the map 
\begin{align*}
	\theta_f: \mathcal{X} \ni x \mapsto \theta_fx\coloneqq \{f_n(x) \}_{n =1}^\infty \in \ell^1(\mathbb{N})
\end{align*}
is a well-defined bounded linear operator and 
\begin{align*}
	|f_n(\tau_n)|\geq1, \quad \forall n \in \mathbb{N}.
\end{align*}
If $x \in \mathcal{X}$ can be written as  $	x=\theta_\tau c$ for some  $c \in \ell^1(\mathbb{N})$ satisfying 
\begin{align*}
	\|c\|_0<\frac{1}{2}\left(1+\frac{1}{\displaystyle\sup_{n, m \in \mathbb{N}, n \neq m}| f_n(\tau_m)|}\right),	
\end{align*}
then $c$ is the unique solution to Problem \ref{BP0}.	
\end{theorem}
\begin{proof}
Theorem \ref{TDP} says that  $c$ is the unique solution to Problem \ref{BP1}.  Let $d \in \ell^1(\mathbb{N})$ be such that 	$x=\theta_\tau d$. We claim that $\|d\|_0> \|c\|_0$.  If this fails, we must have $\|d\|_0\leq \|c\|_0$. We then have 
\begin{align*}
	\|d\|_0<\frac{1}{2}\left(1+\frac{1}{\displaystyle\sup_{n, m \in \mathbb{N}, n \neq m}| f_n(\tau_m)|}\right).
\end{align*}
Theorem \ref{TDP} again says  that $d$ is also the unique solution to Problem \ref{BP1}. Therefore we must have  $\|c\|_1<\|d\|_1$ and $\|c\|_1>\|d\|_1$ which is a contradiction. Therefore the claim holds and we have $\|d\|_0> \|c\|_0$.
\end{proof}
\begin{corollary}
Theorem \ref{DEGN} follows from Theorems \ref{TDP}  and 	\ref{M}.
\end{corollary}
\begin{proof}
 	Let $\{\tau_j\}_{j=1}^n$ be a  normalized frame  for  a Hilbert space $\mathcal{H}$. For each $1\leq j \leq n$, define 
 	\begin{align*}
 		f_j:\mathcal{H}\ni h \mapsto f_j(h)\coloneqq \langle h, \tau_j \rangle \in \mathbb{K}.
 	\end{align*}
 Then 
 \begin{align*}
 	|f_j(\tau_j)|=1, \quad \forall 1\leq j \leq n
 \end{align*}
and 
\begin{align*}
	f_k(\tau_j)=\langle \tau_j, \tau_k\rangle, \quad \forall 1\leq j, k \leq n.
\end{align*}
\end{proof}

\section{Acknowledgments}
The author thanks the anonymous reviewer for his/her suggestions, which resulted in a significant generalization of the main result. The simplified title of this article is also due to the reviewer.

 \bibliographystyle{plain}
 \bibliography{reference.bib}

\end{document}